\documentclass [12pt]{article}
\usepackage{graphicx,amssymb,amsfonts,latexsym,amsmath,amsthm,times}
\usepackage{epsfig}
\usepackage{color}
\usepackage[a4paper]{geometry}
\def\lanbox{\hbox{$\, \vrule height 0.25cm width 0.25cm depth 0.01cm \,$}}

\numberwithin{equation}{section}

\begin{document}

\vspace*{1.4cm}

\normalsize \centerline{\Large \bf SOLVABILITY
OF A CLASS OF INTEGRO-}

\medskip

\centerline{\Large\bf DIFFERENTIAL EQUATIONS WITH  LAPLACE}

\medskip

\centerline{\Large\bf   AND BI-LAPLACE OPERATORS}

\vspace*{1cm}

\centerline{\bf Vitali Vougalter$^{1 *}$,  Vitaly Volpert$^{2,3}$,  }

\vspace*{0.5cm}

\centerline{$^{1 *}$ Department of Mathematics, University
of Toronto}

\centerline{Toronto, Ontario, M5S 2E4, Canada}

\centerline{ e-mail: vitali@math.toronto.edu}

\medskip

\centerline{$^2$ Institute Camille Jordan, UMR 5208 CNRS,
University Lyon 1}

\centerline{ Villeurbanne, 69622, France}

\centerline{$^3$ Peoples' Friendship University of Russia, 6 Miklukho-Maklaya
St,}

\centerline{Moscow, 117198, Russia}

\centerline{e-mail: volpert@math.univ-lyon1.fr}

\medskip


\vspace*{0.25cm}

\noindent {\bf Abstract:}
The work deals with the studies of the existence of solutions of an
integro-differential equation in the situation of the difference of the standard Laplacian
and the bi-Laplacian in the diffusion term. The proof
of the existence of solutions relies on a fixed point technique. We use the solvability
conditions for the non-Fredholm elliptic operators in unbounded domains.

\vspace*{0.25cm}



\vspace*{0.5cm}

\bigskip

\bigskip


\setcounter{equation}{0}

\section{\bf Introduction}

\medskip

The present article is devoted to the studies of the existence of the stationary
solutions of the following integro-differential equation
\begin{equation}
\label{h}
\frac{\partial u}{\partial t} =
D[\Delta-{\Delta}^{2}]u +
\int_{{\mathbb R}^{d}}K(x-y)g(u(y,t))dy + f(x), \quad 5\leq d\leq 7. 
\end{equation}
The problems of this kind are crucial to the cell population dynamics. The results of the work
are obtained in this particular range of the values of the dimension
which is based on the solvability of the linear Poisson type equation (\ref{lp}) and the
applicability of the Sobolev embedding (\ref{e}). 
Existence of stationary solutions in $H^{3}({\mathbb R}^{5})$ for the equation analogous to (\ref{h}) containing
a single Laplace operator in the diffusion term was covered in
~\cite{VV15} (see also ~\cite{VV111}).
The space variable $x$ here
corresponds to the cell genotype,
$u(x,t)$ stands for the cell density as a function of their genotype and time.

The right side of our equation describes the evolution of the cell density by
means of the cell proliferation, mutations and cell influx/efflux.
The diffusion term involving the difference of the Laplace operator and the bi-Laplacian in our context corresponds to
the change of the genotype due to the small random mutations, and the integral
production term describes large
mutations. The function $g(u)$ denotes the rate of the cell birth depending
on $u$
(density dependent proliferation), and the kernel $K(x-y)$ designates
the proportion of the newly born cells changing their genotype from
$y$ to $x$.
It is assumed here that it depends on the distance between the genotypes.
The last term in the right side of (\ref{h}) stands for
the influx or efflux of cells for different genotypes.

The bi-Laplacian in the diffusion term is relevant to the long range interactions in biological systems (see ~\cite{M03}).
The diffusion operator containing the difference of the Laplacian and the bi-Laplacian can arise in several modelling contexts. In evolutionary dynamics, the Laplacian term describes small random mutations in the genotype space, while the higher-order bi-Laplacian may represent long-range smoothing effects or higher-order mutation processes. Similar operators appear in ecological models where individuals disperse through a combination of local diffusion and longer-range redistribution. From another perspective, the operator $\Delta - \Delta^2$
 can be viewed as a local approximation of nonlocal dispersal operators obtained by Taylor expansion of convolution kernels. Finally, higher-order diffusion terms frequently arise in pattern formation theory and in models with higher-gradient regularization, where they reflect additional spatial interactions or energetic penalties associated with curvature of the density distribution.

We set here $D=1$ and demonstrate the existence of solutions of the problem
\begin{equation}
\label{p}
[\Delta-{\Delta}^{2}]u +\int_{{\mathbb R}^{d}} K(x-y)g(u(y))dy +
f(x) = 0, \quad 5\leq d\leq 7.
\end{equation}
Let us discuss the case when the linear part of this operator does not
satisfy the Fredholm property. Consequently,  the conventional methods of
the nonlinear analysis may not be applicable. We use the solvability conditions
for the non-Fredholm operators along with the method of contraction mappings.
Existence of solutions for certain non-Fredholm integro-differential equations with 
the bi-Laplacian was established in ~\cite{VV21}. In ~\cite{EV25} the authors prove
the global well-posedness of an integro-differential equation with the bi-Laplacian
and transport.

Consider the problem
\begin{equation}
\label{eq1}
 -\Delta u + V(x) u - a u=f
\end{equation}
with $u \in E= H^{2}({\mathbb R}^{d})$ and  $f \in F=
L^{2}({\mathbb R}^{d}), \ d\in {\mathbb N}$, $a$ is a constant and
the scalar potential function $V(x)$ is either trivial
or converges to $0$ at infinity. For $a \geq 0$, the essential spectrum of the
operator $A : E \to F$, which corresponds to the left  side of equation
(\ref{eq1}) contains the origin. As a consequence, such  operator fails to
satisfy the Fredholm property. Its image is not closed, for $d>1$
the dimension of its kernel and the codimension of its image are
not finite. The present article deals with the studies of some properties
of the operators of this kind. Note that the elliptic equations with
non-Fredholm operators were studied actively in recent years.
Approaches in weighted Sobolev and H\"older spaces were developed in
~\cite{Amrouche1997}, ~\cite{Amrouche2008}, ~\cite{Bolley1993},
~\cite{Bolley2001}, ~\cite{B88}. The non-Fredholm Schr\"odinger type operators
were treated with the methods of the spectral and the
scattering theory in ~\cite{EV21}, ~\cite{EV22}, ~\cite{V2011}, ~\cite{VV08}, ~\cite{VV14}, ~\cite{VV22}.
The nonlinear non-Fredholm elliptic problems were discussed in ~\cite{EV22}, ~\cite{VV111},
~\cite{VV14}, ~\cite{VV15}, ~\cite{VV21}. Article ~\cite{DH22} deals with a two phase boundary obstacle-type problem for the bi-Laplacian.
In ~\cite{DKV19} the authors studied the limit behaviour of a singular perturbation problem for the biharmonic operator.
The significant applications to the theory of the reaction-diffusion
equations were developed in ~\cite{DMV05}, ~\cite{DMV08}. Fredholm structures,
topological invariants and applications were considered in ~\cite{E09}.
The works ~\cite{GS05} and ~\cite{RS01} are important for the understanding
of the Fredholm and properness properties of the quasilinear elliptic systems
of the second order and of the operators of this kind on ${\mathbb R}^{N}$.
The operators without
the Fredholm property appear also when studying the wave systems with an infinite
number of localized traveling waves (see ~\cite{AMP14}). The standing lattice
solitons in the discrete NLS equation with saturation were covered in
~\cite{AKLP19}.
In particular,
when $a=0$ the operator $A$ is Fredholm in certain properly chosen weighted
spaces (see \cite{Amrouche1997}, \cite{Amrouche2008}, \cite{Bolley1993},
\cite{Bolley2001}, \cite{B88}). However, the situation of $a \neq 0$ is considerably
different and the approach developed in these articles cannot be applied. 

\medskip

Let us set $K(x) = \varepsilon {\cal K}(x)$ with $\varepsilon \geq 0$
and suppose that the conditions below are satisfied.

\medskip

\noindent
{\bf Assumption 1.1.}  {\it Let $f(x): {\mathbb R}^{d}\to {\mathbb R}, \ 5\leq d\leq 7$ be nontrivial, such that
$f(x)\in L^{1}({\mathbb R}^{d})\cap L^{2}({\mathbb R}^{d})$. Assume also that
${\cal K}(x): {\mathbb R}^{d}\to {\mathbb R}$ does not vanish identically in the whole space and
${\cal K}(x)\in L^{1}({\mathbb R}^{d})\cap L^{2}({\mathbb R}^{d})$.}

\bigskip

We choose the space dimension $5\leq d\leq 7$, which is relevant to the solvability
conditions for the linear Poisson type equation (\ref{lp}) formulated in
Lemma 4.1 further down and to the applicability of the Sobolev embedding (\ref{e}). From the point of view of the applications, the space
dimensions are not limited to $5\leq d\leq 7$ because our space variable corresponds to
the cell genotype but not to the usual physical space. 

For the technical purposes, we use the Sobolev space
$$
H^{4}({\mathbb R}^{d}):=\big\{u(x):{\mathbb R}^{d}\to {\mathbb R} \ | \
u(x)\in L^{2}({\mathbb R}^{d}), \ {\Delta}^{2} u \in L^{2}
({\mathbb R}^{d}) \big\}.
$$
It is equipped with the norm
\begin{equation}
\label{n}
\|u\|_{H^{4}({\mathbb R}^{d})}^{2}:=\|u\|_{L^{2}({\mathbb R}^{d})}^{2}+
\big\|{\Delta}^{2}u\big\|_{L^{2}({\mathbb R}^{d})}^{2}.
\end{equation}
By virtue of the standard Sobolev embedding in dimensions $d\leq 7$, we have
\begin{equation}
\label{e}
\|u\|_{L^{\infty}({\mathbb R}^{d})}\leq c_{e}\|u\|_{H^{4}({\mathbb R}^{d})},
\end{equation}
where $c_{e}>0$ is the constant of the embedding. 
When the nonnegative parameter $\varepsilon$ vanishes, we arrive at the
linear Poisson type equation (\ref{lp}) for $x\in {\mathbb R}^{d}, \ 5\leq d\leq 7$.
By means of the result of Lemma 4.1 below along with Assumption 1.1, problem (\ref{lp})
admits a unique solution
\begin{equation}
\label{u0}
u_{0}(x)\in H^{4}({\mathbb R}^{d}), \quad 5\leq d\leq 7.
\end{equation}
Note that (\ref{u0}) does not vanish identically in the whole space since the function $f(x)$ is nontrivial
as assumed.

Let us seek the resulting solution of nonlinear problem (\ref{p}) as
\begin{equation}
\label{r}
u(x)=u_{0}(x)+u_{p}(x).
\end{equation}
Evidently, we derive the perturbative equation
\begin{equation}
\label{pert}
[-\Delta+{\Delta}^{2}]u_{p}(x)=\varepsilon \int_{{\mathbb R}^{d}}
{\cal K}(x-y)g(u_{0}(y)+u_{p}(y))dy, \quad 5\leq d\leq 7. 
\end{equation}
We introduce a closed ball in the Sobolev space, namely
\begin{equation}
\label{b}
B_{\rho}:=\{u(x)\in H^{4}({\mathbb R}^{d}) \ | \ \|u\|_{H^{4}({\mathbb R}^{d})}\leq
\rho \}, \quad 0<\rho\leq 1, \quad 5\leq d\leq 7. 
\end{equation}
Let us look for the solution of equation (\ref{pert}) as the fixed point of the
auxiliary nonlinear problem
\begin{equation}
\label{aux}
[-\Delta+{\Delta}^{2}]u(x)=\varepsilon \int_{{\mathbb R}^{d}}
{\cal K}(x-y)g(u_{0}(y)+v(y))dy, \quad 5\leq d\leq 7 
\end{equation}
in ball (\ref{b}). For a given function $v(y)$ this is an equation with
respect to $u(x)$.
The left side of (\ref{aux}) invloves the operator which fails to satisfy
the Fredholm property
\begin{equation}
\label{l}  
l:=-\Delta+{\Delta}^{2}:
H^{4}({\mathbb R}^{d})\to L^{2}({\mathbb R}^{d}).
\end{equation}
This is the differential operator with the symbol $|p|^{2}+|p|^{4}$, such that
$$
lu(x)=\frac{1}{(2\pi)^{\frac{d}{2}}}\int_{{\mathbb R}^{d}}(|p|^{2}+|p|^{4})
\widehat{u}(p)e^{ipx}dp, \quad u(x)\in H^{4}({\mathbb R}^{d}),
$$
where the standard Fourier transform is introduced in (\ref{f}).
Clearly, the essential spectrum of (\ref{l}) fills the nonnegative semi-axis
$[0, +\infty)$.
Hence, this operator has no bounded inverse. The similar situation
appeared in earlier works ~\cite{VV111} and ~\cite{VV14}. But as distinct from the
present case, the problems considered there required the orthogonality conditions.
Persistence of pulses for certain local reaction-diffusion equations using the fixed 
point technique was covered in ~\cite{CV21}.
But the Schr\"odinger type operator involved in the nonlinear
problem there had the Fredholm property.

Let us use the interval on the real line
\begin{equation}
\label{i}
I:=\big[-c_{e}\|u_{0}\|_{H^{4}({\mathbb R}^{d})}-c_{e}, \
c_{e}\|u_{0}\|_{H^{4}({\mathbb R}^{d})}+c_{e}\big], \quad 5\leq d\leq 7
\end{equation}
along with the closed ball in the space of $C_{2}(I)$ functions, namely
\begin{equation}
\label{M}
D_{M}:=\{g(z)\in C_{2}(I) \ | \ \|g\|_{C_{2}(I)}\leq M \}, \quad M>0.
\end{equation}
The norm involved in (\ref{M})
\begin{equation}
\label{gn}
\|g\|_{C_{2}(I)}:=\|g\|_{C(I)}+\|g'\|_{C(I)}+\|g''\|_{C(I)}
\end{equation}
with $\|g\|_{C(I)}:=\hbox{max}_{z\in I}|g(z)|$.

We impose the following technical conditions on the nonlinear part of
equation  (\ref{p}). It will be trivial at the origin along with its first
derivative. From the perspective of the biological applications, $g(z)$ can
be, for example the quadratic function to describe the cell-cell interaction.

\bigskip

\noindent
{\bf Assumption 1.2.} {\it Let $g(z): {\mathbb R}\to {\mathbb R}$, such that
$g(0)=0$ and $g'(0)=0$. Additionally, we assume that $g(z)\in D_{M}$ and
it does not vanish identically on the interval $I$}.

\bigskip

Let us introduce the operator $t_g$, so that $u = t_g v$, where $u$ is a
solution of problem (\ref{aux}). Our first main result is as follows.

\bigskip

\noindent
{\bf Theorem 1.3.} {\it Let Assumptions 1.1 and 1.2 be valid. Then 
for every $\rho\in (0, 1]$ equation (\ref{aux})
defines the map $t_{g}: B_{\rho}\to B_{\rho}$, which is a strict contraction
for all
$$
0<\varepsilon\leq
$$
\begin{equation}
\label{eps}  
\frac{\rho}{2M(\|u_{0}\|_{H^{4}({\mathbb R}^{d})}+1)^{2}
\Big[\frac{\|{\cal K}\|_{L^{1}({\mathbb R}^{d})}^{2}
(\|u_{0}\|_{H^{4}({\mathbb R}^{d})}+1)^{\frac{8}{d}-2}d}{(2\pi)^{4}(d-4)}
\Big(\frac{|S^{d}|}{16}\Big)^{\frac{4}{d}}
+\frac{\|{\cal K}\|_{L^{2}({\mathbb R}^{d})}^{2}}{4}\Big]^{\frac{1}{2}}}.
\end{equation}
The unique fixed point $u_{p}(x)$ of this map $t_{g}$ is the only solution of
problem (\ref{pert}) in $B_{\rho}$.}

\bigskip

Here and further down $S^{d}$ stands for the unit sphere in the space of $5\leq d \leq 7$
dimensions centered at the origin and $|S^{d}|$ denotes its Lebesgue
measure.

Clearly, the resulting solution of equation (\ref{p}) given by formula (\ref{r})
will be nontrivial because the source term $f(x)$ does not vanish identically in the whole space and
$g(0)=0$ as we assume.

We have the following trivial proposition.

\bigskip

\noindent
{\bf Lemma 1.4.} {\it For $R\in (0, +\infty)$  and $5\leq d\leq 7$ consider the
function
$$
\varphi(R):=\alpha R^{d-4}+\frac{1}{R^{4}}, \quad \alpha>0.
$$
It attains the minimal value at \
$\displaystyle{R^{*}:=\Bigg(\frac{4}{\alpha (d-4)}\Bigg)^{\frac{1}{d}}}$,
which is equal to}
$$
\varphi(R^{*})=\Bigg(\frac{\alpha}{4}\Bigg)^{\frac{4}{d}}\frac{d}
{(d-4)^{\frac{d-4}{d}}}.
$$

\bigskip

Our second main statement deals with the continuity of the resulting solution
of equation (\ref{p}) given by (\ref{r}) with respect to
the nonlinear function $g$.

\bigskip

\noindent
{\bf Theorem 1.5.} {\it Let $j=1,2$, the assumptions of Theorem 1.3 hold,
such that $u_{p,j}(x)$ is the unique fixed point of the map
$t_{g_{j}}: B_{\rho}\to B_{\rho}$, which is a strict contraction for all the values
of $\varepsilon$ satisfying bound (\ref{eps})
and the resulting solution of problem (\ref{p}) with
$g(z)=g_{j}(z)$ is given by
\begin{equation}
\label{ressol}  
u_{j}(x)=u_{0}(x)+u_{p,j}(x).
\end{equation}
Then for all the values of $\varepsilon$, which satisfy inequality (\ref{eps}),
the estimate
$$
\|u_{1}-u_{2}\|_{H^{4}({\mathbb R}^{d})}\leq \frac{\varepsilon}{1-\varepsilon \sigma}
(\|u_{0}\|_{H^{4}({\mathbb R}^{d})}+1)^{2}\times
$$
\begin{equation}
\label{cont}
\Bigg[\frac{\|{\cal K}\|_{L^{1}({\mathbb R}^{d})}^{2}
(\|u_{0}\|_{H^{4}({\mathbb R}^{d})}+1)^{\frac{8}{d}-2}|S^{d}|^{\frac{4}{d}}}
{{16}^{\frac{4}{d}}(2\pi)^{4}}\frac{d}{d-4}+
\frac{\|{\cal K}\|_{L^{2}({\mathbb R}^{d})}^{2}}{4}\Bigg]^{\frac{1}{2}}\|g_{1}-g_{2}\|_{C_{2}(I)}
\end{equation}
is valid.}

\bigskip

We have $\sigma$ introduced in formula (\ref{sig}) further down.
Let us proceed to the proof of our first main proposition.

\bigskip


\setcounter{equation}{0}

\section{\bf The existence of the perturbed solution}

\bigskip

\noindent
{\it Proof of Theorem 1.3.} Let us choose arbitrarily $v(x)\in B_{\rho}$ and
designate the term involved in the integral expression in the right side of
problem (\ref{aux}) as
$$
G(x):=g(u_{0}(x)+v(x)).
$$
We use the standard Fourier transform
\begin{equation}
\label{f}
\widehat{\phi}(p):=\frac{1}{(2\pi)^{\frac{d}{2}}}\int_{{\mathbb R}^{d}}\phi(x)
e^{-ipx}dx, \quad 5\leq d\leq 7.
\end{equation}
Obviously, the estimate from above
\begin{equation}
\label{fub}
\|\widehat{\phi}(p)\|_{L^{\infty}({\mathbb R}^{d})}\leq \frac{1}{(2\pi)^{\frac{d}{2}}}
\|\phi(x)\|_{L^{1}({\mathbb R}^{d})}
\end{equation}
is valid. Let us apply (\ref{f}) to both sides of equation (\ref{aux}) to obtain
$$
\widehat{u}(p)=\varepsilon (2\pi)^{\frac{d}{2}}
\frac{\widehat{\cal K}(p)\widehat{G}(p)}{|p|^{2}+|p|^{4}}
$$
Hence, for the norm we derive
\begin{equation}
\label{un}
\|u\|_{L^{2}({\mathbb R}^{d})}^{2}={(2\pi)}^{d} \varepsilon^{2}\int_{{\mathbb R}^{d}}
\frac{|\widehat{\cal K}(p)|^{2}|\widehat{G}(p)|^{2}}
{[|p|^{2}+|p|^{4}]^{2}}dp.
\end{equation}
Let us use the analog of inequality (\ref{fub}) applied to functions ${\cal K}$ and $G$ with
$R>0$, such that
$$
{(2\pi)}^{d} \varepsilon^{2}\int_{{\mathbb R}^{d}}
\frac{|\widehat{\cal K}(p)|^{2}|\widehat{G}(p)|^{2}}
{[|p|^{2}+|p|^{4}]^{2}}dp\leq 
$$
$$
{(2\pi)}^{d} \varepsilon^{2}\int_{|p|\leq R}\frac{|\widehat{\cal K}(p)|^{2}|
\widehat{G}(p)|^{2}}{|p|^{4}}dp+
{(2\pi)}^{d} \varepsilon^{2}\int_{|p|>R}\frac{|\widehat{\cal K}(p)|^{2}
|\widehat{G}(p)|^{2}}{|p|^{4}}dp\leq
$$
\begin{equation}
\label{ub1}
\varepsilon^{2}\|{\cal K}\|_{L^{1}({\mathbb R}^{d})}^{2}\Bigg\{
\frac{1}{(2\pi)^{d}}\|G(x)\|_{L^{1}({\mathbb R}^{d})}^{2}|S^{d}|\frac{R^{d-4}}
{d-4}+
\frac{1}{R^{4}}\|G(x)\|_{L^{2}({\mathbb R}^{d})}^{2}\Bigg\}.
\end{equation}
Note that for $v(x)\in B_{\rho}$ we have
$$
\|u_{0}+v\|_{L^{2}({\mathbb R}^{d})}\leq \|u_{0}\|_{H^{4}({\mathbb R}^{d})}+1.
$$
Sobolev embedding (\ref{e}) implies that
$$
|u_{0}+v|\leq c_{e}(\|u_{0}\|_{H^{4}({\mathbb R}^{d})}+1).
$$
Evidently,
$$
G(x)=\int_{0}^{u_{0}+v}g'(z)dz.
$$
Thus,
$$
|G(x)|\leq \hbox{sup}_{z\in I}|g'(z)||u_{0}+v|\leq M|u_{0}+v|
$$
with the interval $I$ introduced in (\ref{i}).
This means that 
$$
\|G(x)\|_{L^{2}({\mathbb R}^{d})}\leq M\|u_{0}+v\|_{L^{2}({\mathbb R}^{d})}\leq M
(\|u_{0}\|_{H^{4}({\mathbb R}^{d})}+1).
$$
Clearly,
$$
G(x)=\int_{0}^{u_{0}+v}dy\Big[\int_{0}^{y}g''(z)dz \Big].
$$
Then
$$
|G(x)|\leq \frac{1}{2}\hbox{sup}_{z\in I}|g''(z)||u_{0}+v|^{2}\leq
\frac{M}{2}|u_{0}+v|^{2}.
$$
Therefore,
\begin{equation}
\label{G1}
\|G(x)\|_{L^{1}({\mathbb R}^{d})}\leq \frac{M}{2}\|u_{0}+v\|_{L^{2}({\mathbb R}^{d})}^{2}
\leq \frac{M}{2}(\|u_{0}\|_{H^{4}({\mathbb R}^{d})}+1)^{2}.
\end{equation}
Hence, we obtain the upper bound for the right side of (\ref{ub1}) equal to
$$
{\varepsilon}^{2}\|{\cal K}\|_{L^{1}({\mathbb R}^{d})}^{2}M^{2}
(\|u_{0}\|_{H^{4}({\mathbb R}^{d})}+1)^{2}
\Bigg\{\frac{(\|u_{0}\|_{H^{4}({\mathbb R}^{d})}+1)^{2}|S^{d}|R^{d-4}}{4 (2\pi)^{d}
(d-4)}+\frac{1}{R^{4}}\Bigg\}
$$
with $R\in (0, +\infty)$. Lemma 1.4 gives us the minimal value
of the expression above, so that 
\begin{equation}
\label{ul2ub}
\|u\|_{L^{2}({\mathbb R}^{d})}^{2}\leq
 {\varepsilon}^{2}\|{\cal K}\|
_{L^{1}({\mathbb R}^{d})}^{2}M^{2}(\|u_{0}\|_{H^{4}({\mathbb R}^{d})}+1)^{2+\frac{8}{d}}
\Bigg(\frac{|S^{d}|}{16}\Bigg)^{\frac{4}{d}}\frac{d}{(2\pi)^{4}(d-4)}.
\end{equation}
Let us recall (\ref{aux}).  Thus,
$$
[-\Delta+{\Delta}^{2}]u(x)=\varepsilon
\int_{{\mathbb R}^{d}}{\cal K}(x-y)G(y)dy, \quad 5\leq d\leq 7.
$$
We use the standard Fourier transform (\ref{f}), the analog of estimate
(\ref{fub}) applied to function $G$ and inequality (\ref{G1}) to arrive at
\begin{equation}
\label{32}
\|{\Delta}^{2} u\|_{L^{2}({\mathbb R}^{d})}^{2}\leq \varepsilon^{2}\|G\|_
{L^{1}({\mathbb R}^{d})}^{2}\|{\cal K}\|_{L^{2}({\mathbb R}^{d})}^{2} \leq \varepsilon^{2}\frac{M^{2}}{4}
(\|u_{0}\|_{H^{4}({\mathbb R}^{d})}+1)^{4}\|{\cal K}\|_{L^{2}({\mathbb R}^{d})}^{2}.
\end{equation}
By virtue of the definition of the norm (\ref{n}) 
along with bounds (\ref{ul2ub}) and (\ref{32}), we obtain that
$$
\|u\|_{H^{4}({\mathbb R}^{d})}\leq \varepsilon (\|u_{0}\|_{H^{4}({\mathbb R}^{d})}+1)^{2}M
\times
$$
\begin{equation}
\label{uh3e}  
 \Bigg[\frac{\|{\cal K}\|_{L^{1}({\mathbb R}^{d})}^{2}(\|u_{0}\|_{H^{4}
({\mathbb R}^{d})}+1)^{\frac{8}{d}-2}d}{(2\pi)^{4}(d-4)} 
\Bigg(\frac{|S^{d}|}{16}\Bigg)^{\frac{4}{d}}
+\frac{ \|{\cal K}\|_{L^{2}({\mathbb R}^{d})}^{2}}{4}\Bigg]^{\frac{1}{2}}\leq \rho
\end{equation}
for all the values of the parameter $\varepsilon$, which satisfy inequality (\ref{eps}).
This means that $u(x)\in B_{\rho}$ as well.

Suppose that for some $v(x)\in B_{\rho}$ there exist
two solutions $u_{1,2}(x)\in B_{\rho}$ of problem (\ref{aux}). Evidently, the
difference function $w(x):=u_{1}(x)-u_{2}(x) \in H^{4}({\mathbb R}^{d})$ solves
the homogeneous equation
$$
[-\Delta+{\Delta}^{2}]w=0.
$$
Note that the operator $l: H^{4}({\mathbb R}^{d})\to L^{2}({\mathbb R}^{d})$
defined in (\ref{l}) does not have any 
nontrivial zero modes. Then $w(x)$ will vanish in
${\mathbb R}^{d}$. Therefore, problem (\ref{aux}) defines a map
$t_{g}: B_{\rho}\to B_{\rho}$ for all the values of $\varepsilon$
satisfying bound (\ref{eps}).

Our goal is to demonstrate that under the given conditions such map is a strict
contraction. We choose arbitrarily $v_{1,2}(x)\in B_{\rho}$. The reasoning above
yields that
$u_{1,2}:=t_{g}v_{1,2}\in B_{\rho}$ as well for $\varepsilon$, which satisfy inequality
(\ref{eps}). Equation (\ref{aux}) gives us
\begin{equation}
\label{aux1}
[-\Delta+{\Delta}^{2}]u_{1}(x)=\varepsilon \int_{{\mathbb R}^{d}}
{\cal K}(x-y)g(u_{0}(y)+v_{1}(y))dy,
\end{equation}
\begin{equation}
\label{aux2}
[-\Delta+{\Delta}^{2}]u_{2}(x)=\varepsilon \int_{{\mathbb R}^{d}}
{\cal K}(x-y)g(u_{0}(y)+v_{2}(y))dy
\end{equation}
with $5\leq d\leq 7$. Define
$$
G_{1}(x):=g(u_{0}(x)+v_{1}(x)), \quad G_{2}(x):=g(u_{0}(x)+v_{2}(x)).
$$
Apply the standard Fourier transform (\ref{f}) to both sides of
problems (\ref{aux1}) and (\ref{aux2}). This yields
$$
\widehat{u_{1}}(p)=\varepsilon (2\pi)^{\frac{d}{2}}
\frac{\widehat{\cal K}(p)\widehat{G_{1}}(p)}{|p|^{2}+|p|^{4}}, \quad
\widehat{u_{2}}(p)=\varepsilon (2\pi)^{\frac{d}{2}}
\frac{\widehat{\cal K}(p)\widehat{G_{2}}(p)}{|p|^{2}+|p|^{4}}.
$$
Then
\begin{equation}
\label{u12n2d}  
\|u_{1}-u_{2}\|_{L^{2}({\mathbb R}^{d})}^{2}=\varepsilon^{2}{(2\pi)}^{d}
\int_{{\mathbb R}^{d}}\frac{|\widehat{\cal K}(p)|^{2}
|{\widehat{G_{1}}(p)}-{\widehat{G_{2}}(p)}|^{2}}{[|p|^{2}+|p|^{4}]^{2}}dp.
\end{equation}
The right side of (\ref{u12n2d}) can be easily estimated from above by means of
inequality (\ref{fub}) as
$$
\varepsilon^{2}{(2\pi)}^{d}\int_{|p|\leq R}\frac{|\widehat{\cal K}(p)|^{2}
|{\widehat{G_{1}}(p)}-{\widehat{G_{2}}(p)}|^{2}}{|p|^{4}}dp+
$$
$$  
\varepsilon^{2}{(2\pi)}^{d}\int_{|p|>R}\frac{|\widehat{\cal K}(p)|^{2}
|{\widehat{G_{1}}(p)}-{\widehat{G_{2}}(p)}|^{2}}{|p|^{4}}dp \leq 
$$
$$
\varepsilon^{2}\|{\cal K}\|_{L^{1}({\mathbb R}^{d})}^{2}\Bigg\{\frac{|S^{d}|}
{(2\pi)^{d}}\|G_{1}(x)-G_{2}(x)\|_{L^{1}({\mathbb R}^{d})}^{2}\frac{R^{d-4}}{d-4}+
\frac{\|G_{1}(x)-G_{2}(x)\|_{L^{2}({\mathbb R}^{d})}^{2}}{R^{4}}\Bigg\}
$$
with $R\in (0,+\infty)$. Obviously,
$$
G_{1}(x)-G_{2}(x)=\int_{u_{0}+v_{2}}^{u_{0}+v_{1}}g'(z)dz.
$$
Hence,
$$
|G_{1}(x)-G_{2}(x)|\leq \hbox{sup}_{z\in I}|g'(z)||v_{1}(x)-v_{2}(x)|\leq
M|v_{1}(x)-v_{2}(x)|.
$$
Therefore,
$$
\|G_{1}(x)-G_{2}(x)\|_{L^{2}({\mathbb R}^{d})}\leq M\|v_{1}-v_{2}\|_
{L^{2}({\mathbb R}^{d})}\leq M\|v_{1}-v_{2}\|_{H^{4}({\mathbb R}^{d})}.
$$
Let us use the equality
$$
G_{1}(x)-G_{2}(x)=\int_{u_{0}+v_{2}}^{u_{0}+v_{1}}dy \Big[\int_{0}^{y}g''(z)dz \Big].
$$
Evidently, $G_{1}(x)-G_{2}(x)$ can be trivially bounded from above in the absolute
value by
$$
\frac{1}{2}\hbox{sup}_{z\in I}|g''(z)||(v_{1}-v_{2})(2u_{0}+
v_{1}+v_{2})|\leq\frac{M}{2}|(v_{1}-v_{2})(2u_{0}+v_{1}+v_{2})|.
$$
By virtue of the Schwarz inequality, we derive the estimate from above for the norm
$$
\|G_{1}(x)-G_{2}(x)\|_{L^{1}({\mathbb R}^{d})}\leq
\frac{M}{2}\|v_{1}-v_{2}\|_{L^{2}({\mathbb R}^{d})}\|2u_{0}+v_{1}+v_{2}\|_
{L^{2}({\mathbb R}^{d})}\leq
$$
\begin{equation}
\label{g12}
M\|v_{1}-v_{2}\|_{H^{4}({\mathbb R}^{d})}
(\|u_{0}\|_{H^{4}({\mathbb R}^{d})}+1).
\end{equation}
Thus, we obtain the upper bound for
$\|u_{1}(x)-u_{2}(x)\|_{L^{2}({\mathbb R}^{d})}^{2}$ equal to
$$
\varepsilon^{2}\|{\cal K}\|
_{L^{1}({\mathbb R}^{d})}^{2}M^{2}\|v_{1}-v_{2}\|_{H^{4}({\mathbb R}^{d})}^{2}
\Big\{\frac{|S^{d}|}{(2\pi)^{d}}(\|u_{0}\|_{H^{4}({\mathbb R}^{d})}+1)^{2}
\frac{R^{d-4}}{d-4}+\frac{1}{R^{4}}\Big\}.
$$
Let us use Lemma 1.4 to minimize such quantity over
$R\in (0,+\infty)$. We arrive at
$$
\|u_{1}(x)-u_{2}(x)\|_{L^{2}({\mathbb R}^{d})}^{2}\leq
$$
\begin{equation}
\label{u12n}
\varepsilon^{2}\|{\cal K}\|_{L^{1}({\mathbb R}^{d})}^{2}M^{2}
\|v_{1}-v_{2}\|_{H^{4}({\mathbb R}^{d})}^{2}\frac{|S^{d}|^{\frac{4}{d}}
(\|u_{0}\|_{H^{4}({\mathbb R}^{d})}+1)^{\frac{8}{d}}}
{(2\pi)^{4}{4}^{\frac{4}{d}}}\frac{d}{d-4}.
\end{equation}
By means of formulas (\ref{aux1}) and (\ref{aux2}), we have
$$
[-\Delta+{\Delta}^{2}](u_{1}(x)-u_{2}(x))=\varepsilon
\int_{{\mathbb R}^{d}}{\cal K}(x-y)[G_{1}(y)-G_{2}(y)]dy.
$$
Let us apply the standard Fourier transform (\ref{f}) along with estimates
(\ref{fub}) and (\ref{g12}). This yields
$$
\|{\Delta}^{2}(u_{1}-u_{2})\|_{L^{2}({\mathbb R}^{d})}^{2}\leq \varepsilon^{2}\|{\cal K}\|_{L^{2}({\mathbb R}^{d})}^{2}
\|G_{1}-G_{2}\|_{L^{1}({\mathbb R}^{d})}^{2}\leq
$$
\begin{equation}
\label{d12}
\varepsilon^{2}\|{\cal K}\|_{L^{2}({\mathbb R}^{d})}^{2}M^{2}
\|v_{1}-v_{2}\|_{H^{4}({\mathbb R}^{d})}^{2}(\|u_{0}\|_{H^{4}({\mathbb R}^{d})}+1)^{2}.
\end{equation}
According to (\ref{u12n}) and (\ref{d12}), the norm
$\|u_{1}-u_{2}\|_{H^{4}({\mathbb R}^{d})}$ can be bounded from above by the expression
$$
\varepsilon M(\|u_{0}\|_{H^{4}({\mathbb R}^{d})}+1)
\Bigg\{\frac{\|{\cal K}\|_{L^{1}({\mathbb R}^{d})}^{2}|S^{d}|^{\frac{4}{d}}
(\|u_{0}\|_{H^{4}({\mathbb R}^{d})}+1)^{\frac{8}{d}-2}}
{(2\pi)^{4}{4}^{\frac{4}{d}}}\frac{d}{d-4}+\|{\cal K}\|_{L^{2}({\mathbb R}^{d})}^{2}\Bigg\}^{\frac{1}{2}}
\times
$$
\begin{equation}
\label{contr}           
\|v_{1}-v_{2}\|_{H^{4}({\mathbb R}^{d})}.
\end{equation}
It can be easily deduced from condition (\ref{eps}) that the constant in 
(\ref{contr}) is less than one.
Therefore, the map $t_{g}: B_{\rho}\to B_{\rho}$ defined by problem
(\ref{aux}) is a strict contraction for all the values of $\varepsilon$, which
satisfy (\ref{eps}). Its unique fixed point $u_{p}(x)$ is the
only solution of equation (\ref{pert}) in the ball $B_{\rho}$. Let us recall
(\ref{uh3e}). Hence, $\|u_{p}(x)\|_{H^{4}({\mathbb R}^{d})}\to 0$ as
$\varepsilon\to 0$.
The resulting
$u(x)\in H^{4}({\mathbb R}^{d})$ introduced in (\ref{r}) solves
problem (\ref{p}). \hfill\lanbox

\bigskip

We turn our attention to demonstrating the validity of the second main
statement of the work.

\bigskip


\setcounter{equation}{0}

\section{\bf The continuity of the resulting solution}

\bigskip

\noindent
{\it Proof of Theorem 1.5.} Clearly, for all the values of the parameter
$\varepsilon$ satisfying condition (\ref{eps}),  we have
$$
u_{p,1}=t_{g_{1}}u_{p,1}, \quad u_{p,2}=t_{g_{2}}u_{p,2}.
$$
Thus,
$$
u_{p,1}-u_{p,2}=t_{g_{1}}u_{p,1}-t_{g_{1}}u_{p,2}+t_{g_{1}}u_{p,2}-
t_{g_{2}}u_{p,2}.
$$
Evidently,
$$
\|u_{p,1}-u_{p,2}\|_{H^{4}({\mathbb R}^{d})}\leq\|t_{g_{1}}u_{p,1}-t_{g_{1}}
u_{p,2}\|_{H^{4}({\mathbb R}^{d})}+\|t_{g_{1}}u_{p,2}-t_{g_{2}}u_{p,2}\|_
{H^{4}({\mathbb R}^{d})}.
$$
According to estimate (\ref{contr}), 
$$
\|t_{g_{1}}u_{p,1}-t_{g_{1}}u_{p,2}\|_{H^{4}({\mathbb R}^{d})}\leq \varepsilon
\sigma\|u_{p,1}-u_{p,2}\|_{H^{4}({\mathbb R}^{d})}.
$$
We have $\varepsilon \sigma<1$ since the map
$t_{g_{1}}: B_{\rho}\to  B_{\rho}$ is a strict contraction under our
assumptions. Here the positive constant
$$
\sigma:=M(\|u_{0}\|_{H^{4}({\mathbb R}^{d})}+1)\times
$$
\begin{equation}
\label{sig}  
\Bigg\{\frac{\|{\cal K}\|_{L^{1}({\mathbb R}^{d})}^{2}|S^{d}|^{\frac{4}{d}}(\|u_{0}\|
_{H^{4}({\mathbb R}^{d})}+1)^{\frac{8}{d}-2}}
{(2\pi)^{4}{4}^{\frac{4}{d}}}\frac{d}{d-4}+\|{\cal K}\|_{L^{2}({\mathbb R}^{d})}^{2}\Bigg\}^{\frac{1}{2}}.
\end{equation}
Hence,
\begin{equation}
\label{sigma}
(1-\varepsilon \sigma)\|u_{p,1}-u_{p,2}\|_{H^{4}({\mathbb R}^{d})}\leq
\|t_{g_{1}}u_{p,2}-t_{g_{2}}u_{p,2}\|_{H^{4}({\mathbb R}^{d})}.
\end{equation}
Let us make use of the fact that $t_{g_{2}}u_{p,2}=u_{p,2}$. Consider
$\gamma(x):=t_{g_{1}}u_{p,2}$. Then
\begin{equation}
\label{12}
[-\Delta+{\Delta}^{2}]\gamma(x)=\varepsilon \int_{{\mathbb R}^{d}}
{\cal K}(x-y)g_{1}(u_{0}(y)+u_{p,2}(y))dy,
\end{equation}
\begin{equation}
\label{22}
[-\Delta+{\Delta}^{2}]u_{p,2}(x)=\varepsilon \int_{{\mathbb R}^{d}}
{\cal K}(x-y)g_{2}(u_{0}(y)+u_{p,2}(y))dy
\end{equation}
with
$5\leq d \leq 7$. We define
$$
G_{1, 2}(x):=g_{1}(u_{0}(x)+u_{p,2}(x)), \quad G_{2, 2}(x):=g_{2}(u_{0}(x)+u_{p,2}(x)).
$$
Let us apply the standard Fourier transform (\ref{f}) to both sides of
problems (\ref{12}) and (\ref{22}). This implies
$$
\widehat{\gamma}(p)=\varepsilon (2 \pi)^{\frac{d}{2}} \frac{\widehat{\cal K}(p)
\widehat{G_{1, 2}}(p)}{|p|^{2}+|p|^{4}}, \quad
\widehat{u_{p,2}}(p)=\varepsilon (2 \pi)^{\frac{d}{2}} \frac{\widehat{\cal K}(p)
\widehat{G_{2, 2}}(p)}{|p|^{2}+|p|^{4}},
$$
such that
\begin{equation}
\label{xiup2n}  
\|\gamma(x)-u_{p,2}(x)\|_{L^{2}({\mathbb R}^{d})}^{2}=\varepsilon^{2}{(2\pi)}^{d}
\int_{{\mathbb R}^{d}}\frac{|{\widehat{\cal K}}(p)|^{2}|\widehat{G_{1, 2}}(p)-
\widehat{G_{2, 2}}(p)|^{2}}{[|p|^{2}+|p|^{4}]^{2}}dp.
\end{equation}
We use  inequality (\ref{fub}) to derive the upper bound on the right side of (\ref{xiup2n}) as
$$
\varepsilon^{2}{(2\pi)}^{d}
\int_{|p|\leq R}\frac{|{\widehat{\cal K}}(p)|^{2}|\widehat{G_{1, 2}}(p)-
\widehat{G_{2, 2}}(p)|^{2}}{|p|^{4}}dp+
$$
$$
\varepsilon^{2}{(2\pi)}^{d}
\int_{|p|>R}\frac{|{\widehat{\cal K}}(p)|^{2}|\widehat{G_{1, 2}}(p)-
\widehat{G_{2, 2}}(p)|^{2}}{|p|^{4}}dp\leq 
$$
$$
\varepsilon^{2}\|{\cal K}\|_{L^{1}({\mathbb R}^{d})}^{2}\Bigg\{\frac{|S^{d}|}
{(2\pi)^{d}}\|G_{1, 2}-G_{2, 2}\|_{L^{1}({\mathbb R}^{d})}^{2}\frac{R^{d-4}}{d-4}+
\frac{\|G_{1, 2}-G_{2, 2}\|_{L^{2}({\mathbb R}^{d})}^{2}}{R^{4}}\Bigg\},
$$
where $R\in (0, +\infty)$. Obviously, the identity
$$
G_{1, 2}(x)-G_{2, 2}(x)=\int_{0}^{u_{0}(x)+u_{p,2}(x)}[g_{1}'(z)-g_{2}'(z)]dz
$$
is valid, so that
$$
|G_{1, 2}(x)-G_{2, 2}(x)|\leq \hbox{sup}_{z\in I}|g_{1}'(z)-g_{2}'(z)|
|u_{0}(x)+u_{p,2}(x)|\leq
$$
$$
\|g_{1}-g_{2}\|_{C_{2}(I)}|u_{0}(x)+u_{p,2}(x)|.
$$
Therefore,
$$
\|G_{1, 2}-G_{2, 2}\|_{L^{2}({\mathbb R}^{d})}\leq \|g_{1}-g_{2}\|_{C_{2}(I)}
\|u_{0}+u_{p,2}\|_{L^{2}({\mathbb R}^{d})}\leq
$$
$$
\|g_{1}-g_{2}\|_{C_{2}(I)}(\|u_{0}\|_{H^{4}({\mathbb R}^{d})}+1).
$$
Let us make use of another equality
$$
G_{1, 2}(x)-G_{2, 2}(x)=\int_{0}^{u_{0}(x)+u_{p,2}(x)}dy\Big[\int_{0}^{y}
(g_{1}''(z)-g_{2}''(z))dz\Big].
$$
Clearly,
$$
|G_{1, 2}(x)-G_{2, 2}(x)|\leq \frac{1}{2}\hbox{sup}_{z\in I}|g_{1}''(z)-g_{2}''(z)|
|u_{0}(x)+u_{p,2}(x)|^{2}\leq
$$
$$
\frac{1}{2}\|g_{1}-g_{2}\|_{C_{2}(I)}|u_{0}(x)+u_{p,2}(x)|^{2}.
$$
Then
$$
\|G_{1, 2}-G_{2, 2}\|_{L^{1}({\mathbb R}^{d})}\leq \frac{1}{2}
\|g_{1}-g_{2}\|_{C_{2}(I)}\|u_{0}+u_{p,2}\|_{L^{2}({\mathbb R}^{d})}^{2}\leq
$$
\begin{equation}
\label{G1222}
\frac{1}{2}\|g_{1}-g_{2}\|_{C_{2}(I)}(\|u_{0}\|_{H^{4}({\mathbb R}^{d})}+1)^{2}.
\end{equation}
This yields the estimate from above for the norm
$\|\gamma(x)-u_{p,2}(x)\|_{L^{2}({\mathbb R}^{d})}^{2}$ given by
$$
\varepsilon^{2}\|{\cal K}\|_{L^{1}({\mathbb R}^{d})}^{2}
(\|u_{0}\|_{H^{4}({\mathbb R}^{d})}+1)^{2}\|g_{1}-g_{2}\|_{C_{2}(I)}^{2}\times
$$
\begin{equation}
\label{xiup2n2}  
\Big[\frac{|S^{d}|(\|u_{0}\|_{H^{4}({\mathbb R}^{d})}+1)^{2}}{4(2\pi)^{d}}
\frac{R^{d-4}}{d-4}+\frac{1}{R^{4}}\Big].
\end{equation}
Let us minimize expression (\ref{xiup2n2}) over $R\in (0, +\infty)$
using Lemma 1.4. Hence,
$\displaystyle{\|\gamma(x)-u_{p,2}(x)\|_{L^{2}({\mathbb R}^{d})}^{2}\leq}$
$$
\varepsilon^{2}
\|{\cal K}\|_{L^{1}({\mathbb R}^{d})}^{2}(\|u_{0}\|_{H^{4}({\mathbb R}^{d})}+1)^{2+\frac{8}{d}}
\|g_{1}-g_{2}\|_{C_{2}(I)}^{2}\frac{|S^{d}|^{\frac{4}{d}}}{16^{\frac{4}{d}}
(2\pi)^{4}}\frac{d}{d-4}.
$$
We recall formulas (\ref{12}) and (\ref{22}), namely
$$
[-\Delta+{\Delta}^{2}]\gamma(x)=\varepsilon \int_{{\mathbb R}^{d}}
{\cal K}(x-y)G_{1, 2}(y)dy,
$$
$$
[-\Delta+{\Delta}^{2}]u_{p,2}(x)=\varepsilon \int_{{\mathbb R}^{d}}
{\cal K}(x-y)G_{2, 2}(y)dy,
$$
so that
$$
[-\Delta+{\Delta}^{2}][\gamma(x)-u_{p,2}(x)]=\varepsilon \int_{{\mathbb R}^{d}}{\cal K}(x-y)
[G_{1, 2}(y)-G_{2, 2}(y)]dy.
$$
Let us apply the standard Fourier transform (\ref{f}) along with bounds
(\ref{fub}) and (\ref{G1222}). This gives us
$$
\|{\Delta}^{2}[\gamma(x)-u_{p,2}(x)]\|_{L^{2}({\mathbb R}^{d})}^{2}\leq
$$
$$
\varepsilon^{2}\|G_{1,2}-G_{2,2}\|_{L^{1}({\mathbb R}^{d})}^{2}\|{\cal K}\|_{L^{2}({\mathbb R}^{d})}^{2}  \leq
\frac{\varepsilon^{2}\|{\cal K}\|_{L^{2}({\mathbb R}^{d})}^{2}}{4}(\|u_{0}\|_{H^{4}({\mathbb R}^{d})}+1)^{4}
\|g_{1}-g_{2}\|_{C_{2}(I)}^{2}.
$$
Thus,
$$
\|\gamma(x)-u_{p,2}(x)\|_{H^{4}({\mathbb R}^{d})}\leq \varepsilon
\|g_{1}-g_{2}\|_{C_{2}(I)}(\|u_{0}\|_{H^{4}({\mathbb R}^{d})}+1)^{2}\times
$$
$$
\Bigg[
\frac{\|{\cal K}\|_{L^{1}({\mathbb R}^{d})}^{2}(\|u_{0}\|_{H^{4}
({\mathbb R}^{d})}+1)^{\frac{8}{d}-2}|S^{d}|^{\frac{4}{d}}}
{{16}^{\frac{4}{d}}(2\pi)^{4}}\frac{d}{d-4}+\frac{\|{\cal K}\|_{L^{2}({\mathbb R}^{d})}^{2}}{4}\Bigg]^
\frac{1}{2}.
$$
We use inequality (\ref{sigma}) to derive that
$$
\|u_{p,1}-u_{p,2}\|_{H^{4}({\mathbb R}^{d})}\leq
\frac{\varepsilon}{1-\varepsilon \sigma}(\|u_{0}\|_{H^{4}({\mathbb R}^{d})}+1)^{2}\times
$$
$$
\Bigg[\frac{\|{\cal K}\|_{L^{1}({\mathbb R}^{d})}^{2}
(\|u_{0}\|_{H^{4}({\mathbb R}^{d})}+1)^{\frac{8}{d}-2}|S^{d}|^{\frac{4}{d}}}
{{16}^{\frac{4}{d}}(2\pi)^{4}}\frac{d}{d-4}+
\frac{\|{\cal K}\|_{L^{2}({\mathbb R}^{d})}^{2}}{4}\Bigg]^{\frac{1}{2}}\|g_{1}-g_{2}\|_{C_{2}(I)}.
$$
Let us complete the proof of our theorem by virtue of formula (\ref{ressol}).
\hfill\lanbox

\bigskip


\setcounter{equation}{0}

\section{Auxiliary results}

\bigskip
\bigskip

We establish the solvability for the linear
Poisson type equation with the sum of the negative Laplacian and the bi-Laplacian in the left side and a square integrable right side
\begin{equation}
\label{lp}
[-\Delta+\Delta^{2}]u=f(x), \quad x\in {\mathbb R}^{d}, \quad d\geq 5.
\end{equation}
The technical proposition below can be easily proved by
applying the standard Fourier transform (\ref{f}) to both sides of 
(\ref{lp}).

\bigskip

\noindent
{\bf Lemma 4.1.} {\it  Let  $f(x): {\mathbb R}^{d}\to
{\mathbb R}, \ d\geq 5$ be nontrivial and
$f(x)\in L^{1} ({\mathbb R}^{d})\cap L^{2} ({\mathbb R}^{d})$.
Then problem (\ref{lp}) admits a unique solution
$u(x)\in H^{4}({\mathbb R}^{d})$.}

\bigskip

\noindent
{\it Proof.} It can be trivially checked that if
$u(x)\in L^{2} ({\mathbb R}^{d})$ solves equation (\ref{lp})
with a square integrable right side, it will be contained
$H^{4}({\mathbb R}^{d})$ as well. To establish that, we apply the standard
Fourier transform (\ref{f}) to both sides of  problem (\ref{lp}). This gives us
$$
(|p|^{2}+|p|^{4})\widehat{u}(p)=\widehat{f}(p)\in L^{2} ({\mathbb R}^{d}).
$$
Hence,
$$
\int_{{\mathbb R}^{d}}[|p|^{2}+|p|^{4}]^{2}|\widehat{u}(p)|^{2}dp<\infty.
$$
Clearly, the equality
$$
\|{\Delta}^{2}u\|_{L^{2} ({\mathbb R}^{d})}^{2}=\int_{{\mathbb R}^{d}}|p|^{8}
|\widehat{u}(p)|^{2}dp<\infty
$$
holds. This means that
${\Delta}^{2}u\in L^{2} ({\mathbb R}^{d})$. Let us recall the definition
of the norm (\ref{n}). We obtain that $u(x)\in H^{4}({\mathbb R}^{d})$.

Let us demonstrate the uniqueness of solutions for our equation. Suppose that
problem (\ref{lp}) has two solutions
$u_{1, 2}(x)\in H^{4}({\mathbb R}^{d})$. Then their difference
$w(x):=u_{1}(x)-u_{2}(x)\in H^{4}({\mathbb R}^{d})$. It satisfies the homogeneous
equation
$$
[-\Delta+{\Delta}^{2}]w=0.
$$
Since the operator
$l: H^{4}({\mathbb R}^{d})\to L^{2}({\mathbb R}^{d})$ introduced in (\ref{l})
does not possess
any nontrivial zero modes, the function $w(x)$ will vanish in ${\mathbb R}^{d}$.

By applying the standard Fourier transform (\ref{f}) to both sides of equation
(\ref{lp}), we arrive at
\begin{equation}
\label{lpf}  
\widehat{u}(p)=\frac{\widehat{f}(p)}{|p|^{2}+|p|^{4}}
\chi_{\{|p|\leq 1\}}+
\frac{\widehat{f}(p)}{|p|^{2}+|p|^{4}}
\chi_{\{|p|>1\}}.
\end{equation}
Here and further down $\chi_{A}$ will denote the characteristic function of a set
$A\subseteq {\mathbb R}^{d}$. Note that the second term in the right side of
(\ref{lpf}) can be estimated from above in the absolute value by
$\displaystyle{\frac{|{\widehat f}(p)|}{2}\in L^{2}({\mathbb R}^{d})}$ due to the
one of the given conditions.

Let us recall inequality (\ref{fub}).
Thus, the first term in the right side of (\ref{lpf}) can be bounded
from above in the absolute value by
\begin{equation}
\label{fps2}
\frac{\|f(x)\|_{L^{1}({\mathbb R}^{d})}}{(2\pi)^{\frac{d}{2}}|p|^{2}}\chi_{\{|p|\leq 1\}}.  
\end{equation}
It can be easily verified that expression
(\ref{fps2}) belongs to $L^{2}({\mathbb R}^{d})$ for $d\geq 5$.
\hfill\lanbox

\bigskip

Let us introduce the correspoding sequence of the approximate equations
related to problem (\ref{lp}) with $n\in {\mathbb N}$, namely
\begin{equation}
\label{lpn}
[-\Delta+{\Delta}^{2}]u_{n}=f_{n}(x), \quad
x\in {\mathbb R}^{d}, \quad d\geq 5.
\end{equation}
The right sides of (\ref{lpn}) tend to the right side of (\ref{lp}) in the appropriate function spaces as
$n\to \infty$. We demonstrate that under the reasonable technical assumptions
each equation (\ref{lpn}) has a unique solution
$u_{n}(x)\in H^{4}({\mathbb R}^{d})$, limiting problem (\ref{lp}) admits
a unique solution $u(x)\in H^{4}({\mathbb R}^{d})$ and
$u_{n}(x)\to u(x)$ in $H^{4}({\mathbb R}^{d})$ as $n\to \infty$. This is
the so called {\it solvability in the sense of sequences} for equation
(\ref{lp}) (see ~\cite{V2011},  also  ~\cite{EV21}, ~\cite{EV22}, ~\cite{VV22}).  The final statement of the work is as follows.

\bigskip

\noindent
{\bf Lemma 4.2.} {\it  Let  $n\in {\mathbb N}, \  f_{n}(x): {\mathbb R}^{d}\to {\mathbb R}, \ d\geq 5$ and
$f_{n}(x)\in L^{1} ({\mathbb R}^{d})\cap L^{2} ({\mathbb R}^{d})$, so that
$f_{n}(x)\to f(x)$ in $L^{1} ({\mathbb R}^{d})$ and
$f_{n}(x)\to f(x)$ in $L^{2} ({\mathbb R}^{d})$ as $n\to \infty$ and $f(x)$ is nontrivial.
Then problems (\ref{lp}) and (\ref{lpn}) possess unique solutions
$u(x)\in H^{4}({\mathbb R}^{d})$ and $u_{n}(x)\in H^{4}({\mathbb R}^{d})$
respectively, such that $u_{n}(x)\to u(x)$ in $H^{4}({\mathbb R}^{d})$ as
$n\to \infty$.}

\bigskip

\noindent
{\it Proof.} According to the result of Lemma 4.1 above, equations
(\ref{lp}) and (\ref{lpn}) have unique solutions
$u(x)\in H^{4}({\mathbb R}^{d})$ and
$u_{n}(x)\in H^{4}({\mathbb R}^{d}), \ n\in {\mathbb N}$ respectively.

Suppose that $u_{n}(x)\to u(x)$ in $L^{2}({\mathbb R}^{d})$ as
$n\to \infty$. It can be easily verified that $u_{n}(x)\to u(x)$ in
$H^{4}({\mathbb R}^{d})$ as $n\to \infty$ as well. Obviously, by means of
(\ref{lpn}) and (\ref{lp}) we have
$$
[-\Delta+{\Delta}^{2}](u_{n}(x)-u(x))=f_{n}(x)-f(x).
$$
Using the standard Fourier transform (\ref{f}), we deduce that
$$
\|{\Delta}^{2}(u_{n}(x)-u(x))\|_{L^{2}({\mathbb R}^{d})}\leq
\|f_{n}(x)-f(x)\|_{L^{2}({\mathbb R}^{d})}\to 0, \quad n\to \infty
$$
as assumed. Recall the norm definition (\ref{n}). Hence,
$u_{n}(x)\to u(x)$ in $H^{4}({\mathbb R}^{d})$ as $n\to \infty$.

Let us apply (\ref{f}) to both sides of problems (\ref{lp}) and (\ref{lpn}). This
yields
\begin{equation}
\label{unuhp}  
\widehat{u_{n}}(p)-\widehat{u}(p)=\frac{\widehat{f_{n}}(p)-\widehat{f}(p)}
{|p|^{2}+|p|^{4}}\chi_{\{|p|\leq 1\}}+\frac{\widehat{f_{n}}(p)-\widehat{f}(p)}
{|p|^{2}+|p|^{4}}\chi_{\{|p|>1\}}.        
\end{equation}
Clearly, the second term in the right side of (\ref{unuhp}) can be bounded
from above in the absolute value by
$\displaystyle{\frac{|\widehat{f_{n}}(p)-\widehat{f}(p)|}{2}}$. This means that
$$
\Bigg\|\frac{\widehat{f_{n}}(p)-\widehat{f}(p)}{|p|^{2}+|p|^{4}}
\chi_{\{|p|>1\}}\Bigg\|_{L^{2}({\mathbb R}^{d})}\leq \frac{1}{2}
\|f_{n}(x)-f(x)\|_{L^{2}({\mathbb R}^{d})}\to 0, \quad n\to \infty
$$
as we assume.

Evidently, the first term in the right side of (\ref{unuhp}) can be estimated
from above in the absolute value by virtue of inequality (\ref{fub}) by
$$
\frac{\|f_{n}(x)-f(x)\|_{L^{1}({\mathbb R}^{d})}}{(2\pi)^{\frac{d}{2}}|p|^{2}}
\chi_{\{|p|\leq 1\}}.
$$
Thus,
$$
\Bigg\|\frac{\widehat{f_{n}}(p)-\widehat{f}(p)}{|p|^{2}+|p|^{4}}
\chi_{\{|p|\leq 1\}}\Bigg\|_{L^{2}({\mathbb R}^{d})}\leq
\frac{\|f_{n}(x)-f(x)\|_{L^{1}({\mathbb R}^{d})}}{(2\pi)^{\frac{d}{2}}}
\sqrt{\frac{|S^{d}|}{d-4}}\to 0, \quad n\to \infty
$$
due to the given condition, such that
$$
u_{n}(x)\to u(x) \quad in \quad L^{2}({\mathbb R}^{d}), \quad d\geq 5 \quad as \quad
n\to \infty,
$$
which completes the proof of our lemma.
\hfill\lanbox

\bigskip


\section*{Acknowledgements}

The first author is grateful to Israel Michael Sigal for the partial support
by the NSERC grant NA 7901.
The second author has been supported by the RUDN University Strategic Academic
Leadership Program.

\end{document}